\documentclass{conm-p-l}

\usepackage{amssymb}

\usepackage[hypertex]{hyperref} 

\newtheorem{theorem}{Theorem}[section]
\newtheorem{cor}[theorem]{Corollary}
\newtheorem{prop}[theorem]{Proposition}
\newtheorem{quest}[theorem]{Question}

\newtheorem{obs}[theorem]{Observation}

\theoremstyle{definition}
\newtheorem{definition}[theorem]{Definition}
\newtheorem{probl}[theorem]{Problem}
\newtheorem{prog}[theorem]{Program}
\newtheorem{example}[theorem]{Example}

\theoremstyle{remark}
\newtheorem{remark}[theorem]{Remark}

\numberwithin{equation}{section}

\begin{document}

\title{Special Functions 
 in Minimal Representations}
\dedicatory{Dedicated to Igor Frenkel
 on the occasion of his 60th birthday with great admire}

\author{Toshiyuki Kobayashi}
\address{Kavli IPMU and Graduate School of Mathematical Sciences, 
The University of Tokyo,
3-8-1 Komaba, Meguro, Tokyo, 153-8914, Japan.}
\email{toshi@ms.u-tokyo.ac.jp}
\thanks{The author was partially supported by
        Grant-in-Aid for Scientific Research (B) (22340026), Japan
        Society for the Promotion of Sciences.}


\subjclass[2010]{Primary 22E30; 
Secondary 22E46, 
33C45
}

\date{}

\begin{abstract}
Minimal representations of a real reductive group $G$ are
the \lq{smallest}\rq\ irreducible unitary representations of $G$. 
We discuss special functions
 that arise
 in the analysis
 of $L^2$-model of minimal representations.  
\end{abstract}

\maketitle

\section{Introduction}
\label{label:intro}
\setcounter{subsection}{1}

An irreducible unitary representation 
 of a real reductive Lie group $G$ 
 is called {\textit{minimal}}
 if its infinitesimal representation is annihilated
 by the Joseph ideal 
 \cite{xjoseph}
in the enveloping algebra.  
Loosely speaking,
 minimal representations
 of $G$ are the \lq{smallest}\rq\
 infinite dimensional
 unitary representations.

The {\it{Weil representation}},
known for a prominent role in number theory
 (e.g.\ the theta correspondence),
provides minimal representations
 of the metaplectic group
 $Mp(n,{\mathbb{R}})$.  
The minimal representation
 of a conformal group $SO(4,2)$
 appears
 in mathematical physics, 
 e.g., 
 as the bound states
 of the hydrogen atom,
 and incidentally
 as the quantum Kepler problem.  
In these classical examples 
 the representations are
 highest weight modules,
 however,
 for more general reductive groups,
 minimal representations
 (if exist)
 may not be highest weight modules,
 see a pioneering work
 of Kostant \cite{xKo} for $SO(4,4)$.

In the last decade
 I have been developing a geometric and analytic theory
 of minimal representations
 with my collaborators,
 S. Ben Sa{\"i}d, J. Hilgert, G. Mano, 
 J. M{\"o}llers, B. {\O}rsted, and M. Pevzner, 
 see \cite{xbko, xhkmm1, xhkmm2, xkhmL2, xhkmo, 
xkcheck, xkgolden, xkmano1, 
xkmano2, xkmanoAMS, xkmollers, 
xkors1, KOP}.  
Among all,
 in this paper, 
 we focus on 
 \lq{special functions}\rq\
 that arise naturally 
 in the $L^2$-model
 of minimal representations.  
Needless to say, 
the interaction
 between special functions and group representations
 has a long history
 and there is extensive literature on this subject.
A new feature in our setting for minimal representations is that
 the representation of the group is realized on the Hilbert space
 $L^2(\Xi)$ where the dimension of a manifold $\Xi$ (see \eqref{eqn:Xi},
 or more generally Section \ref{subsec:schrodinger} for the definition of $\Xi$)
 is strictly smaller than the dimension of any nontrivial
 $G$-space in most cases.
This means that
 $G$ cannot act geometrically on $\Xi$ but there is a natural action 
 of $G$ on $L^2(\Xi)$.
As a result, the Casimir element of a compact subgroup
acts as a \textit{fourth-order} differential operator. 
The indefinite orthogonal group $G=O(p+1, q+1)$
 is the most interesting for this purpose 
 in the sense
 that the group $G$ itself contains two parameters
 $p$ and $q$, 
 and we shall highlight 
 this case
 by giving occasionally some perspectives 
 to other reductive groups.

\vskip 0.8pc
\par\noindent
{\bf{Acknowledgements:}}\enspace
I first met Igor Frenkel when I visited Yale University 
 in 2009 to give a colloquium talk.
It was then a surprising pleasure 
 that Igor told me his recent theory
 on quaternionic analysis \cite{xfm08, xfm2011, FL12} with Libine uses some 
of my work \cite{xkmanoAMS, xkors1}
 on geometric analysis of minimal representations, 
which encourages me to develop further the analytic theory of minimal 
representations.

I would like to thank the organizers P. Etingof, M. Khovanov, A. 
Kirillov Jr., A. Lachowska, A. Licata, A. Savage and G. Zuckerman for 
their hospitality 
during the stimulating conference \lq\lq Perspectives in Representation Theory\rq\rq\ in honor of Prof. Igor Frenkel's 60th 
birthday
at Yale University,  12--16 May 2012.
Thanks are also due to an anonymous referee for his/her careful comments.

\section{A generalization of the Fourier transform}
\subsection{Algebraic characterization
 of Fourier transforms}
\label{subsec:algchara}
We begin
 with an algebraic characterization
 of the Euclidean Fourier transform 
 ${\mathcal{F}}_{\mathbb{R}^n}$.  
Let $Q_j:=x_j$
 be the multiplication operators
 by coordinates,
 and $P_j:=\frac{1}{\sqrt{-1}}\frac{\partial}{\partial x_j}$.  
Then we have:
\begin{prop}
\label{prop:F}
Any continuous operator $A$ on $L^2({\mathbb{R}}^n)$
 satisfying
\[
  A \circ Q_j = P_j \circ A, 
  \quad
  A \circ P_j = - Q_j \circ A
  \quad
  \text{ on } {\mathcal{S}}({\mathbb{R}}^n)
  \quad
  (1 \le j \le n)
\]
 is a scalar multiple
 of the Euclidean Fourier transform ${\mathcal {F}}_{\mathbb{R}^n}$.  
In particular,
 any such continuous operator $A$ is unitary up to scaling.  
\end{prop}

In place of ${\mathbb{R}}^n$, 
 let us consider the isotropic cone
\begin{equation}
\label{eqn:Xi}
\Xi:=\{
x \in {\mathbb{R}}^{p+q}\setminus \{0\}:
x_{1}^2 + \cdots + x_{p}^2 - x_{p+1}^2 - \cdots -x_{p+q}^2
=0
\}, 
\end{equation}
equipped 
 with a measure $d \mu = \frac 1 2 r^{p+q-3} d r d \omega d \eta$
 in the bipolar coordinates:
\begin{equation}
\label{eqn:bipolar}
{\mathbb{R}}_+ \times S^{p-1} \times S^{q-1}
\overset \sim \to \Xi,
\quad
 (r,\omega,\eta) \mapsto (r \omega, r \eta).  
\end{equation}

\begin{remark}
This cone $\Xi$ is a special case of the Lagrangian submanifold
 of a minimal real coadjoint orbit,
 denoted by the same letter $\Xi$,
 given in Theorem \ref{thm:khmL2} where we deal with more general reductive groups.
\end{remark}

Recall from \cite{xkmanoAMS}
 that the {\textit{fundamental differential operators}}
 $R_j$ ($1 \le j \le p+q$)
 on $\Xi$ are mutually commuting operators 
 which are obtained as the restriction 
 of the tangential differential operators
\[
  \varepsilon_j x_j \square
 -(2 E + p+ q-2)\frac{\partial}{\partial x_j} 
\]
to $\Xi$, 
where 
$\varepsilon_j =1$ 
 ($1 \le j \le p$); 
 $=-1$
 ($p+1 \le j \le p+q$), 
$
   \square:= \sum_{a=1}^{p+q} \varepsilon_a 
             \frac{\partial^2}{\partial x_a^2}
$ (the Laplacian on ${\mathbb{R}}^{p,q}$)
 and 
$
   E:= \sum_{a=1}^{p+q} x_a 
       \frac{\partial}{\partial x_a}  
$ (the Euler operator).  
Then we have 
\begin{theorem}
[{\cite[Theorem 1.2.3]{xkmanoAMS}}]
\label{thm:AR}
Suppose $p+q$
 is even,
 $\ge 4$.  
Then there exists
 a unitary operator ${\mathcal{F}}_{\Xi}$
 on $L^2(\Xi)$
 satisfying the following relation for $A$:
\begin{equation}
\label{eqn:AR}
  A \circ Q_j = R_j \circ A, 
  \quad
  A \circ R_j = Q_j \circ A
  \quad
  \text{on }\,\, C_c^{\infty}(\Xi)
  \quad
  (1 \le j \le p+q).  
\end{equation}
Conversely,
 any continuous operator $A$
satisfying \eqref{eqn:AR}
 is a scalar multiple
 of the unitary operator ${\mathcal{F}}_{\Xi}$.  
\end{theorem}

\subsection{Unitary inversion operator
 ${\mathcal{F}}_{\Xi}$}
\label{subsec:2.2}

The similar nature
 of ${\mathcal{F}}_{\mathbb{R}^n}$
 and ${\mathcal{F}}_{\Xi}$ 
 indicated in Proposition \ref{prop:F}
 and Theorem \ref{thm:AR}
 is derived from the common fact
 that they arise as the {\textit{unitary inversion operators}}
 in the $L^2$-model
 (Schr{\"o}dinger model)
 of minimal representations
 of real reductive groups
 $Mp(n,{\mathbb{R}})$ and $O(p+1,q+1)$, 
 respectively.  

To see some more details of Theorem \ref{thm:AR}, 
 let $\overline{\mathfrak {n}}$, 
 ${\mathfrak {n}}$, 
 and ${\mathfrak {l}}$ be the Lie algebras
 generated by the operators
 $Q_i$, $R_j$, 
 and $[Q_i,R_j]$
 ($1 \le i,j \le p+q$).  
Then ${\mathfrak {g}}:=\overline{\mathfrak {n}}
 +{\mathfrak {l}}+{\mathfrak {n}}$
 is isomorphic to ${\mathfrak {o}}(p+1,q+1)$, 
 and 
$
   {\mathfrak {p}}
=
   {\mathfrak {l}}
+
   {\mathfrak {n}}
\simeq ({\mathfrak {o}}(p,q)+{\mathbb{R}}) \ltimes {\mathbb{R}}^{p+q}
$
 is a maximal parabolic subalgebra of ${\mathfrak {g}}$.

For $p+q$ even, 
 $\ge 4$, 
we proved in \cite{xkors1}
 that there exists an irreducible unitary representation
 $\pi$ of the group $G:=O(p+1,q+1)$
 on the Hilbert space $L^2(\Xi)$
 of which the infinitesimal representation
 is given by $Q_j$
 (the action of 
$
  {\mathfrak {n}}
$)
and $R_j$ 
(the action of 
$
  {\overline{\mathfrak {n}}}
$), 
see also \cite[Chapter 1]{xkmanoAMS}.  
We set 
\begin{equation}
\label{eqn:Ipq}
   w := \begin{pmatrix} I_p & O \\ O & -I_q
        \end{pmatrix}\in G.  
\end{equation}
Geometrically,
 ${\mathfrak {p}}$ is the Lie algebra of the conformal transformation group
 $(O(p,q)\cdot\mathbb{R}_{>0}) \ltimes\mathbb{R}^{p,q}$ 
 of the flat pseudo-Riemannian 
 Euclidean space 
 ${\mathbb {R}}^{p,q}$, 
 and $w$ induces the conformal inversion of ${\mathbb{R}}^{p,q}$
 by the M{\"o}bius transform.

The unitary operators $\pi(g)$
 are of simple form
 if $g\in G$ belong to the maximal parabolic subgroup $P$
 with Lie algebra ${\mathfrak {p}}$, 
namely,
 they are given by the multiplication
 of certain elementary functions on $\Xi$
 and the translations
 coming from the geometric action
 of the Levi subgroup of $P$ on $\Xi$
 (see \cite[Chap.2, Sect.3]{xkmanoAMS}).  
In view of the Bruhat decomposition 
\[G = P \amalg PwP, 
\]
it is enough to find an explicit formula
 of the unitary operator $\pi(w)$
 in order to give a global formula
 of the $G$-action on $L^2(\Xi)$.  
We call $\pi(w)$ 
 the {\it{unitary inversion operator}}, 
 and set
\begin{equation}\label{eqn:FXi}
\mathcal{F}_\Xi := \pi(w).  
\end{equation}
We initiated
 in a series of papers \cite{xkmano1, xkmano2, xkmanoAMS}
 the following:
\begin{prog}
[{\cite[Program 1.2.5]{xkmanoAMS}}]
\label{prog:FA}
Use the unitary inversion operator
 ${\mathcal {F}}_{\Xi}$
 for minimal representations
 as an analog 
 of the Euclidean Fourier transform ${\mathcal {F}}_{\mathbb{R}^n}$, 
 and develop a theory
 of \lq{Fourier analysis}\rq\
 on $\Xi$.  
\end{prog}

In the classical Schr{\"o}dinger model
 of the Weil representation
 of the metaplectic group $Mp(n,{\mathbb{R}})$
on $L^2({\mathbb{R}}^n)$, 
 the unitary inversion operator
 is nothing but the Euclidean Fourier transform
 ${\mathcal{F}}_{{\mathbb{R}}^n}$
 (up to scalar of modulus one), 
 see Example \ref{ex:F}.  
We note
 that $Mp(n,{\mathbb{R}})$
 and $O(p+1,q+1)$
 with $p+q$ even 
 are simple Lie groups 
 of type $C$ and $D$, 
 respectively.  

\vskip 0.8pc
The first stage of Program \ref{prog:FA}
 is to establish
 a framework of the $L^2$-model
 (Schr{\"o}dinger model)
 of minimal representations, 
and to introduce the unitary inversion operator
 ${\mathcal{F}}_{\Xi}$
 with an algebraic characterization 
 such as Theorem \ref{thm:AR}.  
In \cite{xkhmL2}
we gave such a model
 and defined ${\mathcal{F}}_{\Xi}$
 by using Jordan algebras,
 see Section \ref{subsec:schrodinger}.  
In this case
 $\Xi$ is a Lagrangian subvariety
 of a minimal nilpotent coadjoint orbit
 and the resulting representations
 on $L^2(\Xi)$ include a slightly wider family
 of unitary representations
 than minimal representations
(e.g.\ the full complementary series representations
 of $O(n,1)$).  
\vskip 0.8pc
The second stage is to solve the following:
\begin{probl}
\label{prog:6.1}
Find an explicit formula
 of the integral kernel
 of ${\mathcal{F}}_{\Xi}$.  
\end{probl}

We will discuss Problem \ref{prog:6.1}
 in Section \ref{sec:genfourier}.  
It is noteworthy 
 that I. Frenkel and M. Libine have developed their original theory
 on quaterionic analysis in a series of papers \cite{xfm08,xfm2011,FL12}
 from the viewpoint of representation theory of the conformal group
 $SL(2,\mathbb{H}_{\mathbb{C}})\simeq SL(4,\mathbb{C})$
 and its real forms,
 and have demonstrated a close connection between minimal
 representations of various $O(p,q)$'s
 and quaternionic analysis.
For instance,
 the explicit formula of $\mathcal{F}_\Xi$ for $O(3,3)$
 which was obtained in Kobayashi--Mano \cite{xkmanoAMS}
 is used in \cite{FL12} for the study of the key operator
 ({\textsl{Pl}$_R$ in their notation)
 in the analysis of the space $\mathbb{H}_{\mathbb{R}}$
 of split quaternions.

The third stage and beyond
 will be based on the algebraic property
 (Theorem \ref{thm:AR}) and analytic property
 (Problem \ref{prog:6.1})
 of the unitary inversion operator ${\mathcal{F}}_{\Xi}$.  
Among various, potentially interesting directions
 of the \lq{Fourier analysis}\rq\ on $\Xi$, 
 here are some few topics:  
\begin{enumerate}
\item[$\bullet$]
A theory of holomorphic semigroups 
 was given in Howe \cite{Howe}
 for $Mp(n,{\mathbb{R}})$
 and in Kobayashi--Mano \cite{xkmano2} for $SO(n+1,2)$.  
\item[$\bullet$]
 A deformation theory
 of the Euclidean Fourier transform ${\mathcal{F}}_{\mathbb {R}^n}$
 \cite{xbko}, 
 e.g.\ an interpolation
 between ${\mathcal{F}}_{\mathbb {R}^n}$
 and the unitary inversion operator
 of $SO(n+1,2)$.  
\item[$\bullet$]
A generalization
 of the classical Bargmann--Segal transform.  
See \cite{xhkmo}
 in the case $G/K$ is of tube type.  
\end{enumerate}

\vskip 0.8pc
Stage 1 already includes 
 a solution
 for the Plancherel-type theorem
 of ${\mathcal{F}}_{\Xi}$.  
A natural but open question would be a Paley--Wiener type
 theorem of ${\mathcal{F}}_{\Xi}$:
\begin{quest}
\label{quest:2.5}
Find an explicit characterization 
 of ${\mathcal{F}}_{\Xi}(C_c^{\infty}(\Xi))$.  
\end{quest}

Another important space
 of functions is an analog
 of Schwartz functions.  
For this we may consider:
\begin{definition}
[Schwartz space on $\Xi$]
\label{dfn:Schwartz}
Let ${\mathcal{S}}(\Xi)$ be the Fr{\'e}chet space
 of smooth vectors
 of the unitary representation
 of $G$ on $L^2(\Xi)$.  
\end{definition}
This definition makes sense 
 in a more general setting
 (see Theorem \ref{thm:khmL2}).  
By the general theory
 of unitary representations,
 we have:
\begin{prop}
${\mathcal{F}}_{\Xi}$
 induces automorphisms
 of the Hilbert space $L^2(\Xi)$
 and the Fr{\'e}chet space ${\mathcal{S}}(\Xi)$.  
\begin{alignat*}{4}
{\mathcal{F}}_{\Xi}: & L^2(\Xi) &&\overset \sim \to \,\,&&L^2(\Xi)
\quad&&\text{{\rm{(}}Plancherel type theorem{\rm{)}}}, 
\\
&\,\,\,\, \cup && &&\,\,\,\, \cup&&
\\
&{\mathcal{S}}(\Xi) &&\overset \sim \to \,\,&&{\mathcal{S}}(\Xi)
\quad&&\text{{\rm{(}}Paley--Wiener type theorem{\rm{)}}}.  
\end{alignat*}
\end{prop}
The following question is also open:
\begin{quest}
Find an explicit characterization of ${\mathcal{S}}(\Xi)$.  
\end{quest}

\subsection{Schr\"odinger model of minimal representations}
\label{subsec:schrodinger}

Suppose that $V$ is a real simple Jordan algebra.
Let $G$ and $L$ be the identity components
 of the conformal group
 and the structure group of the Jordan algebra $V$,
 respectively.
Then the Lie algebra ${\mathfrak {g}}$
 is a real simple Lie algebra
 and has a Gelfand--Naimark decomposition
$ {\mathfrak {g}}= \overline{{\mathfrak {n}}}
                    + {\mathfrak {l}}
                    + {\mathfrak {n}},
$
where ${\mathfrak {n}} \simeq V$
 is regarded as an Abelian Lie algebra,
 ${\mathfrak {l}} \simeq {\mathfrak {str}}(V)$
 the structure algebra,
 and $\overline{{\mathfrak {n}}}$
 acts on $V$
 by quadratic vector fields.

Let 
$
   {\mathbb{O}}_{\operatorname{min}}^{G_{\mathbb{R}}}
   (\subset {\mathfrak {g}}^{\ast})
$
 be a (real) minimal nilpotent coadjoint orbit.
By identifying ${\mathfrak {g}}$
 with the dual ${\mathfrak {g}}^{\ast}$,
 we consider the intersection 
$$
\Xi:=V \cap {\mathbb{O}}_{\operatorname{min}}^{G_{\mathbb{R}}},
$$ 
 which is a Lagrangian submanifold
 of the symplectic manifold 
 ${\mathbb{O}}_{\operatorname{min}}^{G_{\mathbb{R}}}$
 endowed with the Kirillov--Kostant--Souriau
 symplectic form.  
There is a natural $L$-invariant Radon measure
 on $\Xi$, 
 and we write $L^2(\Xi)$
 for the Hilbert space 
 consisting of square integrable functions 
 on $\Xi$.  
\begin{theorem}
[Schr\"odinger model \cite{xkhmL2}]
\label{thm:khmL2}
Suppose $V$ is a real simple Jordan algebra such that
 its maximal Euclidean Jordan algebra is also simple.
Among all such Jordan algebras $V$,
 we exclude the case where $V \simeq \mathbb{R}^{p,q}$
 with $p+q$ odd (see Examples \ref{exam:symmsp} and \ref{exam:4.2}).
\par\noindent
{\rm{1)}}\enspace
For an appropriate finite covering group $\widetilde{G}$
 of $G$
 there exists a natural unitary representation 
 of $\widetilde{G}$ on $L^2(\Xi)$.  
It is irreducible
 if and only if 
 $\Xi$ is connected.  
\par\noindent
{\rm{2)}}\enspace
The Gelfand--Kirillov dimension of $\pi$
attains its minimum
among all infinite dimensional representations
 of $\widetilde{G}$, 
 i.e.\ $\operatorname{DIM}(\pi)=\frac 1 2 \dim {\mathbb{O}}_{\operatorname{min}}^{G_{\mathbb{R}}}$.  
\par\noindent
{\rm{3)}}\enspace
The annihilator of the differential representation
 $d \pi$ is the Joseph ideal
 in the enveloping algebra $U({\mathfrak {g}}_{\mathbb{C}})$
 if $V$ is split and ${\mathfrak {g}}_{\mathbb{C}}$
 is not of type $A$.
\end{theorem}

The simple Lie algebras ${\mathfrak {g}}$
 that appear in Theorem \ref{thm:khmL2}
 are categorized into four cases
 as below:
\begin{align}
&{\mathfrak {sl}}(2k,{\mathbb{R}}),
{\mathfrak {so}}(2k,2k),
{\mathfrak {so}}(p+1,q+1),
{\mathfrak {e}}_{7(7)},
\label{eqn:J1}
\\
&
{\mathfrak {sp}}(k,{\mathbb{R}}),
{\mathfrak {su}}(k,k),
{\mathfrak {so}}^{\ast}(4k),
{\mathfrak {so}}(2,k),
{\mathfrak {e}}_{7(-25)},
\label{eqn:J2}
\\
&{\mathfrak {sp}}(k,{\mathbb{C}}),
{\mathfrak {sl}}(2k,{\mathbb{C}}),
{\mathfrak {so}}(4k,{\mathbb{C}}),
{\mathfrak {so}}(k+2,{\mathbb{C}}),
{\mathfrak {e}}_{7}({\mathbb{C}}),
\label{eqn:J3}
\\
&{\mathfrak {sp}}(k,k),
{\mathfrak {su}}^{\ast}(4k),
{\mathfrak {so}}(k,1).
\label{eqn:J4}
\end{align}

\begin{example}
If $V$ is a Euclidean Jordan algebra, 
 then $G$ is the automorphism group
 of a Hermitian symmetric space
 of tube type
 and the corresponding Lie algebra ${\mathfrak {g}}$
 is listed in \eqref{eqn:J2}.  
In this case $\Xi$ consists
 of two connected components, 
 and the resulting representation $\pi$
 is the direct sum
 of an irreducible unitary 
 highest weight module
 and its dual.  
\end{example}
\begin{remark}
In the case \eqref{eqn:J4}
 the complex minimal nilpotent orbit
${\mathbb{O}}_{\operatorname{min}}^{G_{\mathbb{C}}}$
 does not meet the real form ${\mathfrak {g}}$, 
 and there does not exist an admissible representation 
of any Lie group
 with Lie algebra ${\mathfrak {g}}$.  
In particular,
 the representation $\pi$ 
 in Theorem \ref{thm:khmL2}
 is not a minimal representation
 but still one of the \lq{smallest}\rq\
 infinite dimensional representations
 in the sense
 that the Gelfand--Kirillov dimension
 attains its minimum.
\end{remark}

\begin{example}
\label{exam:symmsp}
Let $V= \operatorname{Sym}(n,{\mathbb{R}})$.
Then
${\mathfrak {g}}={\mathfrak{sp}}(n,{\mathbb{R}})$
 and
\begin{equation}
\label{eqn:coneSp}
\Xi=\{X \in M(n,{\mathbb{R}}):
 X = {}^{t\!}X, \operatorname{rank}X=1\}.
\end{equation}
Let 
$
   \Xi_+
:=\{X \in V \cap {\mathbb{O}}_{\operatorname{min}}^{G_{\mathbb{R}}}:
    \operatorname{Trace}X>0\}.
$
Via the double covering map ({\it{folding map}})
\[
  {\mathbb{R}}^n \setminus\{0\} \to \Xi_+,
  \quad
  v \mapsto v {}^{t\!}v
\]
 we can identify the representation on $L^2(\Xi_+)$ with the even part
 of the Schr{\"o}dinger model on $L^2({\mathbb{R}}^n)$
 of the metaplectic group $Mp(n,{\mathbb{R}})$
 \cite{xfolland, Howe}.
See \cite{xkhmL2}
 for the realization
 of the odd part of the Weil representation
 in the space of sections for a certain line bundle over $\Xi_+$.  
\end{example}
\begin{example}
\label{exam:4.2}
We define a multiplication on
$\mathbb{R}^{p+q} = \mathbb{R} \oplus \mathbb{R}^{p+q-1}$
by
\[
(x_1,x') \cdot (y_1,y') :=
   (x_1y_1 - \sum_{i=2}^p x_i y_i + \sum_{i=p+1}^{p+q} x_i y_i,
   x_1y' + y_1x').
\]
The resulting Jordan algebra is denoted by $\mathbb{R}^{p,q}$
(by a little abuse of notation).
It is a semisimple Jordan algebra of rank two,
and its conformal algebra is $\mathfrak{o}(p+1,q+1)$.
Suppose now
 $V= {\mathbb{R}}^{p,q}$
 with $p+q$ even.
Then
 $\Xi$ in Theorem \ref{thm:khmL2} coincides
 with the isotropic cone given in \eqref{eqn:Xi}.  
For $q=1$, 
 $V$ is an Euclidean Jordan algebra,
 and $\Xi$ consists of two connected components
 according to the sign
of the first coordinate $x_1$,
 i.e.~the past and future cones.  
For $p,q \ge 2$, 
 $V$ is non-Euclidean,
 $\Xi$ is connected,
 and our representation $\pi$ on $L^2(\Xi)$ is 
 the same as the Schr{\"o}dinger model
 of the minimal representation
 of $O(p+1,q+1)$
 constructed in \cite[Part III]{xkors1}, 
 which is neither a highest nor a lowest weight module.
\end{example}

\begin{remark}
There is no minimal representation
 for any group
 with Lie algebra ${\mathfrak {o}}(p+1,q+1)$
 with $p+q$ odd,
 $p$, $q$ $\ge 3$
 (see \cite[Theorem 2.13]{xV}).

\end{remark}

\section{Unitary inversion operator ${\mathcal{F}}_{\Xi}$}
\label{sec:genfourier}

By the Schwartz kernel theorem,
 the unitary inversion operator ${\mathcal{F}}_{\Xi}$
 can be given by a distribution kernel
 $K(x,y) \in {\mathcal{D}}'(\Xi \times \Xi)$, 
 namely, 
\[
   {\mathcal{F}}_{\Xi}u(x)
   =
  \int_{\Xi}K(x,y)u(y)d\mu (y)
  \qquad
  \text{for all }\,\,
  u \in C_c^{\infty}(\Xi).  
\]
Problem \ref{prog:6.1} asks for
 an explicit formula
 of $K(x,y)$.  
In the setting of Theorem \ref{thm:khmL2},
 we can generalize the definition \eqref{eqn:FXi} of $\mathcal{F}_\Xi$
 by taking $w$ to be a lift of the conformal inversion on $V$,
 see \cite{xhkmo}.
So far,
 Problem \ref{prog:6.1} has been solved
 for minimal representations
 in the following two cases:
\par\noindent
Case A.\enspace
 $G$ is the biholomorphic transformation group
 of a Hermitian symmetric space
 of tube type
 (\cite{xhkmo}).  
\par\noindent
Case B. \enspace
 $G=O(p+1,q+1)$
 (Theorem \ref{thm:3.2}).  
\par
Case A includes the following earlier results:
\begin{example}
\label{ex:F}
\begin{enumerate}
\item[1)]
$G=Mp(n,\mathbb{R})$, 
 $\Xi={\mathbb{R}}^n$, 
 $\pi=$ the Weil representation.  
\[
 K(x,y)=\frac{c}{(2 \pi)^{\frac n 2}} e^{-\sqrt{-1}\langle x,y \rangle}
\]
In this case ${\mathcal{F}}_{\Xi}$
 is the Euclidean Fourier transform 
 ${\mathcal{F}}_{{\mathbb{R}}^n}$
 up to a phase factor $c$
 with $|c|=1$,
 see \cite{xfolland}.  
\item[2)]
$G=SO(p+1,2)$, 
 $\Xi$ is the light cone
 ($q=1$ in \eqref{eqn:Xi}),
 $\pi=$ the highest weight representation of the smallest
 Gelfand--Kirillov dimension and of the smallest $K$-type.
\[
 K(x,y)= c \widetilde{J}_{\frac {p-3}{2}}(2 \sqrt{2\langle x,y\rangle})
\]
 where $\widetilde{J}_{\lambda}(t):=(\frac{t}{2})^{-\lambda}J_{\lambda}(t)$
 is a renormalization
 of the J-Bessel function
 (\cite{xkmano1}).  
\end{enumerate}
\end{example}

\subsection{Mellin--Barnes type integral expression}

In \cite{xkmanoAMS}
 we brought an idea of the Radon transform in the analysis of the unitary
 inversion operator $\mathcal{F}_\Xi$ for minimal representations.
Recall
 that the Euclidean Fourier transform 
 ${\mathcal{F}}_{\mathbb{R}^n}$
 can be written as the composition of the one-dimensional Fourier transform
 and the Radon transform
 ({\textit{plane wave decomposition}}).  
We can generalize this decomposition to the unitary inversion operator $\mathcal{F}_\Xi$
 for the minimal representations
 and some small representations
 on $L^2(\Xi)$
 given by Theorem \ref{thm:khmL2},
 namely, 
 there exists a distribution $\Phi(t)$ of one variable such that
 the distribution kernel $K(x,y)$
 of ${\mathcal{F}}_{\Xi}$
 is of the following form:
\begin{equation}
\label{eqn:Radon}
  K(x,y)=\Phi(\langle x,y \rangle ),
\end{equation}
 where $\langle\, , \,\rangle$
 is some (natural) bilinear form
 of the ambient space $V$.
Thus Problem \ref{prog:6.1} reduces to find a formula of $\Phi(t)$.

\begin{example}
In Example \ref{ex:F} (2) we have seen
 $\Phi(t) = c \widetilde{J}_{\frac{p-3}{2}} (2 \sqrt{2t})$
 when $G = O(p+1,2)$.
Therefore,
 ${\mathcal{F}}_{\Xi}$
 reduces to the Hankel transform
 composed
 with a \lq{Radon transform}\rq\
 on $\Xi$
 in this case.  
\end{example}

The formula of $\Phi(t)$ is more involved for
 $G = O(p+1,q+1)$ with $p,q \ge 2$
 as the corresponding minimal representation is not a highest weight module and
 $\Phi(t)$ is not always locally integrable
 (see Question \ref{quest:3.4} below).
An explicit formula of $\Phi(t)$ in this case can be given
 in terms of `Bessel distributions'
 \cite[Theorem 5.1.1]{xkmanoAMS}.  
Here we give an alternative expression
 of $\Phi(t)$,
 namely,
 by using a distribution-valued Mellin--Barnes integral.

We define a distribution of $t$
 with meromorphic parameter $\lambda$ by 
\[
   b(\lambda, t):=\frac{\Gamma(-\lambda)}{\Gamma(\lambda+\frac{p+q}{2}-1)}(2t)_+^{\lambda}.  
\]
Here the Riesz distribution $(2t)_+^\lambda$ is defined
 as a locally integrable function on $\mathbb{R}$ by
\[
(2t)_+^\lambda
  = \begin{cases}   (2t)^\lambda   & t > 0 \\
                         0                   & t \le 0
     \end{cases}
\]
for $\operatorname{Re}\lambda > -1$,
and is extended as a distribution by the meromorphic continuation on
$\lambda\in\mathbb{C}$.
Let $m:=\frac 1 2 (p+q-4)$, 
 and $L_m$ be a contour starting
 at $\gamma -i \infty$, 
 passes the real axis
 between $(-m-1,-m)$
 and ends at $\gamma + i \infty$
 when $\gamma >-1$.  
We define distributions $\Phi^{p,q}(t)$
 by a distribution-valued Mellin--Barnes integral:
\begin{equation*}
\Phi^{p,q}(t)
:=
\begin{cases}
\int_{L_0} b(\lambda, t) d \lambda
\qquad
&(\text{Case A-1}), 
\\
\int_{L_m} b(\lambda, t) d \lambda
\qquad
&(\text{Case B-1}), 
\\
\int_{L_m} (\frac{b(\lambda, t)}{\tan \pi \lambda}
+
\frac{b(\lambda, -t)}{\sin \pi \lambda}) d \lambda
\qquad
&(\text{Case B-2}),
\end{cases}
\end{equation*}
according to the following three cases:
\par\noindent
Case A-1. \enspace $p=1$  or $q=1$, 
\par\noindent
Case B-1. \enspace $p, q>1$ and both odd, 
\par\noindent
Case B-2. \enspace $p,q>1$ and both even.  

Then $\Phi^{p,q}(t)$ is independent of the choice of the contour and
 $\gamma$ under the above mentioned constraints.

Let $\langle \ , \ \rangle$ be the (positive definite) inner product
 on $\mathbb{R}^{p,q}$.
Then we have
\begin{theorem}
[{\cite[\S 6.2]{xkmanoAMS}}]
\label{thm:3.2}
For $G=O(p+1,q+1)$ 
 with $p+q$ even, 
 $\ge 4$, 
 the kernel $K(x,y)$
 of the unitary inversion ${\mathcal{F}}_{\Xi}$
 is given by 
\begin{equation*}
K(x,y)=c_{p,q}\Phi^{p,q}(\langle x, y \rangle)  
\end{equation*}
 for some constant $c_{p,q}$.
\end{theorem}

\subsection{Local integrability of the kernel}
The kernel of the Euclidean Fourier transform
 ${\mathcal{F}}_{\mathbb{R}^n}$
 is given by $e^{-i\langle x, \xi \rangle}$, 
 which is locally integrable.  
We may ask 
 to which extent this analytic feature remains to hold.  
To be more precise,
 let $\Phi(t)$ be the distribution
 on ${\mathbb{R}}$ as in \eqref{eqn:Radon}.  
We ask
\begin{quest}
\label{quest:3.4}
When is $\Phi(t)$ locally integrable?  
\end{quest}
For a Euclidean Jordan algebra $V$
 we proved in \cite{xhkmo} that
 $\Phi(t)$ is locally integrable.  
See \eqref{eqn:J2}
 for the list of the corresponding conformal
 Lie algebras ${\frak{g}}$.

For $G=O(p+1, q+1)$
 with $p+q$ even $>2$, 
 the Mellin--Barnes type integral formula (Theorem \ref{thm:3.2})
 leads to the following proposition
 (see \cite[Theorem 6.2.1]{xkmanoAMS}):

\begin{prop}
\label{prop:3.5}
We have the identities modulo
$L_{\operatorname{loc}}^{1}({\mathbb{R}}, r^{p+q-3}dr)$.  
\begin{equation*}
\Phi^{p,q}(t)
\equiv
\begin{cases}
0
\qquad
&{\text{\rm(Case A-1),}}
\\
c_1 \sum_{l=0}^{m-1} \frac{(-1)^l}{2^l(m-l-1)!}
\delta^{(l)}(t)
\qquad
&{\text{\rm{(Case B-1),}}}
\\
c_2 \sum_{l=0}^{m-1} \frac{l!}{2^l(m-l-1)!}
t^{-l-1}
\qquad
&{\text{\rm{(Case B-2),}}}
\end{cases}
\end{equation*}
for some nonzero constants $c_1$, $c_2$. 
Here $m = \frac12(p+q-4)$. 
\end{prop}

Thus we have a complete answer to Question \ref{quest:3.4} in this case:
\begin{cor}
\label{cor:3.6}
$\Phi(t)$ is locally integrable
 if and only if\/
 ${\mathfrak {g}}={\mathfrak {o}}(p+1,2), 
{\mathfrak {o}}(2,q+1)$ or\/ 
 ${\mathfrak {o}}(3,3) \simeq {\mathfrak {sl}}(4,{\mathbb{R}})$.  
\end{cor}

We note that
 the minimal representations
 for 
 $\mathfrak{g} = \mathfrak{o}(p+1,2)$ or $\mathfrak{o}(2,q+1)$
 are highest (or lowest) weight modules,
 whereas minimal representations
 do not make sense
 for $\mathfrak{g} = \mathfrak{o}(3,3)$
 which is isomorphic to $\mathfrak{sl}(4,\mathbb{R})$. 
(Recall that the Joseph ideal 
is defined 
 when ${\mathfrak {g}}_{\mathbb{C}}$
 is not of type $A$.)

The delicate answer indicated in Corollary \ref{cor:3.6}
is closely related to the regularity 
 of the `Radon transform' on $\Xi$.
To be more precise,
 the Radon transform on $\Xi$ is defined as the integral
 over the codimension-one submanifold
\[
 \langle x,y \rangle = t \quad\text{in $\Xi$}.
\]
 which collapses when $t=0$.
Accordingly,
 the Radon transform
\[
 (\mathcal{R} u)(x,t)
 = \int_\Xi u(y) \delta(\langle x,y \rangle - t) dy
\]
 has a better regularity as $t$ tends to $0$.
(On the other hand,
 the asymptotic behavior as $|t| \to \infty$ is similar to the Euclidean case.)
The singular part 
 of $\Phi^{p,q}(t)$ 
 in Proposition \ref{prop:3.5} fits well
 with the behavior of ${\mathcal{R}}u(x,t)$
 as $t$ tends to 0.  

Question \ref{quest:3.4} is open  
 for minimal representations
 without highest weights
 except for the case $G=O(p+1,q+1)$.  

\section{Fourth order differential equations}

\subsection{Gaussian kernel and minimal $K$-type}
\label{subsec:gaussiankernel}

The Euclidean Fourier transform 
 ${\mathcal{F}}_{\mathbb{R}^n}$
 is of order four,
 and therefore its eigenvalues
 in $L^2({\mathbb{R}}^n)$
 are among $\{\pm 1, \pm \sqrt{-1}\}$. 
An important eigenfunction
 with eigenvalue 1 is the Gaussian kernel
 $e^{-\frac{1}{2}\|x\|^2}$,
 namely, 
\[
  {\mathcal{F}}_{\mathbb{R}^n}(e^{-\frac{1}{2}\|x\|^2})
=e^{-\frac{1}{2}\|x \|^2}.  
\]
Thus the Gaussian kernel $e^{-\frac 1 2 \|x\|^2}$
 is a square integrable function on ${\mathbb{R}}^n$
 satisfying the following property:
\begin{equation}
\label{eqn:FRinv}
{\mathcal{F}}_{\mathbb{R}^n}f =f 
\quad\text{and}\quad
\text{$f$ is $O(n)$-invariant.  }
\end{equation}

Next let $p \ge q \ge 1$, 
$p+q$ even, 
 and we consider the isotropic cone $\Xi$
 in $\mathbb{R}^{p,q}$
 as in \eqref{eqn:Xi}.  
The unitary inversion operator ${\mathcal{F}}_{\Xi}$
 of the minimal representation $\pi$
 of $O(p+1,q+1)$
 is of order two
 because ${\mathcal{F}}_{\Xi}=\pi(w)$
 and $w^2=I_{p+q}$
 (see \eqref{eqn:Ipq}), 
 and therefore
 its eigenvalues in $L^2(\Xi)$
 are either $1$ or $-1$.  
An important eigenfunction of ${\mathcal{F}}_{\Xi}$
 is $\widetilde K_{\frac1 2 (q -2)}(2\|x\|)$
 where 
$
  \widetilde K_{\lambda}(t):=(\frac t 2)^{-\lambda}K_{\lambda}(t)
$ 
is a renormalization
 of the K-Bessel function.  
This is a square integrable function on $\Xi$
 satisfying the following property: 
\begin{equation}
\label{eqn:FXinv}
{\mathcal{F}}_{\Xi}f =\pm f 
\quad\text{and}\quad
\text{$f$ is $O(p) \times O(q)$-invariant.  }
\end{equation}
To be more precise,
\begin{equation*}
{\mathcal{F}}_{\Xi}f
=
\begin{cases}
f 
\qquad
&\text{ if } p-q \equiv 0 \mod 4,
\\
-f
&\text{ if } p-q \equiv 2 \mod 4.  
\end{cases}
\end{equation*}

\begin{example}
\label{ex:hydro}
For $p=3$ and $q=1$, 
 we have 
\[
    {\mathcal{F}}_{\Xi}(e^{-2\|x\|})
   =-e^{-2\|x\|}, 
\]
because  
\begin{equation}
\label{eqn:Khalf}
   \widetilde K_{-\frac 1 2}(t)
  =\frac{\sqrt{\pi}}{2}e^{-t}.
\end{equation}  

The function $e^{-2\|x\|}$ arises as the wave function
 for the hydrogen atom 
 with the lowest energy.   
\end{example}

{}From the view point 
 of representation theory, 
 the Gaussian kernel 
 $e^{-\frac 1 2 \|x\|^2}$ generates
 the minimal $K$-type 
 of the Weil representation of $Mp(n,{\mathbb{R}})$, 
 whereas the function $\widetilde{K}_{\frac 1 2 (q-2)}(2\|x\|)$
 generates
 that of the minimal representation of $O(p+1,q+1)$
 realized in $L^2(\Xi)$.

\subsection{The Mano polynomial}
\label{subsec:Mano}
We recall a classical fact 
 that the Hermite polynomials form
 an orthogonal basis
 for the radial part of the Schr{\"o}dinger model $L^2({\mathbb{R}}^n)$
 of the Weil representation e.g.\ \cite{xfolland,Howe},
 whereas the Laguerre polynomials arise in the minimal representation
 of the conformal group $SO(4,2)$.
The bottom of the series
 correspond to what we have discussed 
 in Section \ref{subsec:gaussiankernel}.  

We notice
 that these two minimal representations
 are quite special,
 namely,
 they are highest weight modules.
However,
 for more general reductive groups,
 minimal representations
 are not always highest weight modules,
 and we need new \lq{orthogonal polynomials}\rq\
 and \lq{special functions}\rq\
 to describe a natural basis
 of functions satisfying
 \eqref{eqn:FXinv}
 or alike.  

For $\mu\in\mathbb{C} \setminus \{-1,-2,-3,\cdots\}$
 and $\ell\in\mathbb{N}$, 
 the {\it{Mano polynomials}}
 $\{M_j^{\mu,\ell}(x)\}_{j\in\mathbb{N}}$
 are defined by 
\begin{equation}
 M_j^{\mu,\ell}(x) :=
 \frac{\Gamma(j+\mu+1)}{j!2^\mu\Gamma(j+\frac{\mu+1}{2})}
\left.\frac{\partial^j}{\partial t^j}\right|_{t=0}
G^{\mu,\ell}(t,x), 
\label{eq:DefPoly}
\end{equation}
where the generating function 
 $G^{\mu, l}(t,x)$
 is given by 
\begin{equation}
    G^{\mu,\ell}(t,x) 
 := \frac{\left(\frac{x}{2}\right)^{2\ell+1}e^{\frac{x}{2}}}
         {(1-t)^{\ell+\frac{\mu+3}{2}}}
\widetilde{I}_{\frac{\mu}{2}}\left(\frac{tx}{2(1-t)}\right)
 \widetilde{K}_{\ell+\frac{1}{2}}\left(\frac{x}{2(1-t)}\right). 
\label{eq:GenFct}
\end{equation}
Here 
$
   \widetilde{I}_\alpha(z):=(\frac{z}{2})^{-\alpha}I_\alpha(z)
$
 and 
$
     \widetilde{K}_\alpha(z)=(\frac{z}{2})^{-\alpha}K_\alpha(z)
$
 denote the renormalized $I$- and $K$-Bessel functions. 
With this normalization,
 the polynomial $M_j^{\mu,\ell}(x)$
 is of the following top term:
$$
     M_j^{\mu,\ell}(x)
  =\frac{(-1)^j}{j!}x^{j+\ell}+\textup{lower order terms}.
$$


\begin{example}
[Special values of the Mano polynomial]
\label{ex:BottomCase}
~~
\begin{enumerate}
\item[\textup{(1)}]
The bottom of the series with $j=0$ is related
 with  the $K$-Bessel functions with half-integer parameter:
\begin{equation*}
   M_0^{\mu,\ell}(x)
=
\pi^{-\frac{1}{2}}z^{2l+1}e^{z}
\widetilde{K}_{\ell+\frac{1}{2}}(z)
\,(=
\sum_{k=0}^\ell{\frac{(2\ell-k)!}{k!(\ell-k)!}x^k}
).
  \end{equation*}
\item[\textup{(2)}]
  The polynomials $M_j^{\mu,\ell}(x)$ for $\ell=0$ reduce to the Laguerre polynomials
  \begin{equation*}
   M_j^{\mu,0}(x) = L_j^\mu(x)
   \,(=
\frac{\Gamma(n+\nu+1)}{n!}
\sum_{k=0}^n{(-1)^k
{n\choose k}\frac{x^k}{\Gamma(k+\nu+1)}}
).
\end{equation*}
\item[\textup{(3)}]
The function $M_{j}^{\mu,\ell}(x)$ is not a polynomial
 when $\ell \not \in {\mathbb{N}}$,
 but it is convenient to include 
 the negative integer case.  
In particular,
 for $\ell = -1$, 
 it follows from \cite[Corollary 5.3]{xhkmm1}
 and \cite[Lemma 3.2]{xhkmm2}
 that 
$M_{j}^{\mu,-1}(x)$ is essentially
 the Laguerre polynomial:  
\[
  x M_{j}^{\mu,-1}(x)=L_{j}^{\mu}(x)
  \qquad
  (j \in {\mathbb{N}}).  
\]
\end{enumerate}
\end{example}

Many of the classical orthogonal polynomials
 are obtained as eigenfunctions of self-adjoint
 differential operators of second-order,
 but the Mano polynomials $M_j^{\mu,\ell}(x)$
 are obtained
 as those of fourth-order.  
Indeed this is a requirement from 
 representation theory
 because the Casimir operator (for a compact subgroup)
 acts as a fourth-order differential operator on $\Xi$.
To see this we may recall that
 the Lie algebra $\overline{{\mathfrak {n}}}$
 acts as a second order differential operator
 (e.g.\ the fundamental differential operator
 $R_j$ in \eqref{eqn:AR}).

We begin with a second order differential operator
 on ${\mathbb{R}}$  
\begin{align*}
\mathcal{R}_{\mu,\ell} :=&
(x\frac{d}{d x}+\mu-2\ell-1-\frac x 2)
(x\frac{d}{d x}+\mu-\frac x 2)
-(\frac x 2)^2, 
\intertext{and introduce a fourth order differential operator}
\mathcal{P}_{\mu,\ell} :=&
 \frac{1}{x^2} \mathcal{R}_{\mu,\ell} \mathcal{R}_{0,\ell}.  
\end{align*}

\begin{prop}
[{\cite{xhkmm2}}]
\label{thm:DiffEq}
{\rm{(1)}}\enspace
{\rm{(Differential equation)}}
The Mano polynomial $M_j^{\mu,\ell}(x)$
 is an eigenfunction
 of ${\mathcal{P}}_{\mu,\ell}$:
\begin{equation}
 \mathcal{P}_{\mu,\ell}u = j(j+\mu+1)u.
\label{eq:DiffEqPoly}
\end{equation}
\par\noindent
{\rm{(2)}}\enspace
{\rm{(Completeness)}}
If $\mu\geq2\ell+1$ is an odd integer,
 then $\{ M_j^{\mu,\ell} \}_{j\in\mathbb{N}}$
 forms an orthogonal basis
 of $L^2(\mathbb{R}_+,x^{\mu-2\ell}e^{-x}d x)$.
\end{prop}
There are two proofs 
 of an explicit formula
 for the $L^2$-norm of $M_{j}^{\mu,\ell}(x)$, 
 see \cite[Theorem 2.4]{xhkmm2}
 and \cite[Corollary 4.1]{xkmollers}.

\subsection{$K$-finite vectors
 for minimal representations}
We discuss a relationship
 between the Mano polynomials
 and minimal representations.
Let $\Xi$ be the isotropic cone in $\mathbb{R}^{p,q}$
 as in \eqref{eqn:Xi} with $p+q \ge 4$ and even.

We recall that the space
 of $k$th spherical harmonics
\[
{\mathcal{H}}^k({\mathbb{R}}^n)
:=
\{\varphi \in C^{\infty}(S^{n-1})
:\Delta_{S^{n-1}}\varphi = -k(k+n-2)\varphi\}
\]
 is a finite dimensional vector space, 
on which the orthogonal group
 $O(n)$ acts irreducibly
 by rotations.  
Then we have an irreducible representation
 of $K:=O(p+1) \times O(q+1)$ on 
\[
   V^j:={\mathcal{H}}^j({\mathbb{R}}^{p+1})
        \otimes
        {\mathcal{H}}^{j+\frac{p-q}{2}}({\mathbb{R}}^{q+1}), 
   \quad
   j=0,1,2,\cdots.  
\]
Let $L^2(\Xi)_K$ be the underlying $({\mathfrak {g}}, K)$-module
 of the minimal representation $L^2(\Xi)$ of $G=O(p+1,q+1)$.  
Then by \cite{xkors1}, 
 $L^2(\Xi)_K$ is isomorphic to the multiplicity-free sum 
 $\oplus_{j=0}^{\infty}V^j$ as $K$-modules.  
We write $L^2(\Xi)_j$ 
for the corresponding $K$-irreducible subspace
 of $L^2(\Xi)$.  

With the notation of Section \ref{subsec:2.2},
 $K \cap P \simeq O(p) \times O(q) \times \mathbb{Z}_2$.
We set $M:=O(p) \times O(q)$.
The representation of the whole group $G$
 (or even the maximal compact subgroup $K$)
 on $L^2(\Xi)$ does not 
 come from the geometric action of $G$ on $\Xi$, 
 but the action of the subgroup $M$
 is given by rotation 
 in the argument
 and the bipolar coordinates \eqref{eqn:bipolar} respect the $M$-action.
Correspondingly,
 it is not straightforward
 to find explicit $K$-finite vectors
 in $L^2(\Xi)$, 
 whereas the $M$-invariant functions
 only depend on the radial parameter 
 $r \in {\mathbb{R}}_+$.  
We identify the space $L^2(\Xi)^M$
 of $M$-invariants in $L^2(\Xi)$
 with $L^2({\mathbb{R}}_+,(1/2) r^{p+q-3} d r)$
 as we saw in Section \ref{subsec:algchara}.

\begin{prop}
[{\cite[Section 8]{xhkmm1}}]
\label{prop:4.4}
Let $G=O(2m,2n)$.  
For every $j \in {\mathbb{N}}$, 
 the subspace $L^2(\Xi)_j \cap L^2(\Xi)^{M}$
 is one-dimensional and given by the Mano polynomials:

\[
  u_j^{m,n}(x):=
  x^{-2n+3}e^{-x}M_{j}^{2m-3,2n-3}(x).  
\]
Here $x=2r$ in the bipolar coordinates \eqref{eqn:bipolar} of $\Xi$.
\end{prop}

Therefore the functions
 $\{ u_j^{m,n}(2r)\}_{j\in\mathbb{N}}$
 give a basis 
 of functions $f$
 satisfying \eqref{eqn:FXinv}.
To be more precise,
 we have
\begin{equation*}
{\mathcal{F}}_{\Xi}u_j^{m,n}
=
\begin{cases}
u_j^{m,n} 
\qquad
&\text{ if } n-m \equiv j \mod 2,
\\
-u_j^{m,n}
&\text{ if } n-m \not\equiv j \mod 2.  
\end{cases}
\end{equation*}

\begin{example}
{\rm{1)}}\enspace
The bottom parameter $j=0$
 explains
 that $\widetilde K_{\frac 1 2 (q-2)}(2r)$
 generates the minimal $K$-type
 (Section \ref{subsec:gaussiankernel}).  
\par\noindent
{\rm{2)}}\enspace
The case $n=1$ recovers the classical fact 
 that the Laguerre polynomials generate
 every $K$-finite vectors
 for the minimal representations
 of $SO(2m,2)$.  
\end{example}
The following remarkable observation
 is 
 a consequence
 of Example \ref{ex:BottomCase} (3)
 and Proposition \ref{prop:4.4}:
\begin{obs}
\label{obs:4.6}
Let $G=O(2m,4)$
 with $m \ge 2$.  
Then any $K$-type of $\pi$
 is given by using Laguerre polynomials.  
\end{obs}

We notice
 that the minimal representations
 of $O(2m,4)$ are not highest weight modules.  
It would be interesting to find a different proof for Observation \ref{obs:4.6},
 possibly in connection
 with other areas of mathematics.

In contrast to Proposition \ref{prop:4.4} for
 $G=O(2m,2n)$,
 $K$-finite vectors in $L^2(\Xi)$ cannot be expressed
 by elementary functions for
 $G=O(2m+1,2n+1)$.
For this
 we introduced in \cite{xhkmm1} a family of real analytic functions
 $\{\Lambda_j^{\mu,\nu}(x)\}_{j \in {\mathbb{N}}}$
 on ${\mathbb{R}}_+$ 
 by the generating function: 
\[
  \sum_{j=0}^{\infty} t^j \Lambda_j^{\mu,\nu}(x)
  =
  \frac{1}{(1-t)^{\frac{\mu+\nu+2}{2}}}
  \widetilde I_{\frac \mu 2}(\frac{tx}{1-t})
  \widetilde K_{\frac \nu 2}(\frac{x}{1-t}).  
\]
Then we have the following Proposition:
\begin{prop}
[{\cite[Corollary 8.2]{xhkmm1}}]
\label{prop:4.7}
Let $G=O(p+1,q+1)$
 with $p+q$ even, 
 $>2$.  
Then 
\[
  L^2(\Xi)_{j} \cap L^2(\xi)^{M}={\mathbb{C}}\Lambda_j^{p-2,q-2}(2r)
\]
for all $j \in {\mathbb{N}}$.  
\end{prop}

\begin{cor}
\label{cor:orthL}
For $\mu \ge \nu \ge -1$, 
 $\mu, \nu \in {\mathbb{Z}}$
 with $\mu \equiv \nu \mod 2$
 and $(\mu,\nu) \ne (-1,-1)$,
$\{\Lambda_{j}^{\mu,\nu}(x)\}_{j \in {\mathbb{N}}}$
 forms an orthogonal basis
 in $L^2({\mathbb{R}}_+, x^{\mu+\nu+1}dx)$.  
\end{cor}

\begin{remark}
The function $\Lambda_j^{\mu,\nu}$ can be expressed by
 elementary functions when
 $\nu\in2\mathbb{Z}+1$:
\[
 \Lambda_j^{\mu,2\ell +1} (x)
 = \frac{2^\mu \Gamma(j+\frac{\mu+1}{2})}
           {\Gamma(j+\mu+1)}
    x^{-2\ell -1} e^{-x} M_j^{\mu,\ell} (2x).
\]
Thus Proposition \ref{prop:4.7} includes Proposition \ref{prop:4.4}
 as a special case.
\end{remark}
\begin{remark}
The indefinite orthogonal group $O(p,q)$ has two parameters $p$ and $q$,
 and the corresponding special functions are expected to be the most general.
In fact,
 for other minimal representations and some small representations
 given in Theorem \ref{thm:khmL2},
 analogous results of Proposition \ref{prop:4.7} remain true
 by a specific choice of the parameters $\mu$ and $\nu$.
\end{remark}

\bibliographystyle{amsplain}

\end{document}